\documentclass[12pt,fleqn]{amsart}

\usepackage{amsthm}
\usepackage{graphicx}
\usepackage{color}
\usepackage{amsmath}
\usepackage{amssymb}
\usepackage{amscd}
\usepackage{latexsym}
\usepackage{graphics}

\theoremstyle{plain}
\newtheorem{theorem}{Theorem}[section]
\newtheorem{lemma}[theorem]{Lemma}
\newtheorem{corollary}[theorem]{Corollary}

\theoremstyle{definition}

\theoremstyle{remark}

\def\I#1{{[\![#1]\!]}}            

\newcommand{\T}[1]{\mathcal{T}_{#1}}

\begin{document}
\title[Meta-Fibonacci Sequences]{Meta-Fibonacci Sequences, Binary Trees, and Extremal Compact Codes}

\author[B. Jackson]{Brad Jackson}
\address{Dept. of Mathematics, San Jose State University, USA}

\author[F. Ruskey]{Frank Ruskey}
\address{Dept. of Computer Science, University of Victoria, CANADA}
\thanks{Research supported in part by NSERC}
\urladdr{http://www.cs.uvic.ca/~ruskey}

\maketitle

\begin{abstract}
We look at a family of meta-Fibonacci sequences which arise in
  studying the number of leaves at the largest level in certain
  infinite sequences of binary trees, restricted compositions of
  an integer, and binary compact codes.
For this family of meta-Fibonacci sequences and two families of
  related sequences we derive ordinary generating functions and
  recurrence relations.
Included in these families of sequences are several well-known
  sequences in the Online Encyclopedia of Integer Sequences (OEIS).
\end{abstract}

\section{Introduction}

In a remarkable paper Emily Norwood studied the number
  of ``compact codes" \cite{Norwood}.
A compact code can be thought of as the sorted sequence of level numbers
  of the leaves of an extended binary tree.
She provided a recurrence relation and table of the number
  of trees classified according to their height and their number of leaves.
We will prove that 
  if the outline of this table is considered as an increasing sequence of
  integers, then one of the ``meta-fibonacci" numbers arises, namely
  the one that satisfies the recurrence relation
\[
a(n) = a(x(n)-a(n-1)) + a(y(n)-a(n-2)),
\]
  with $x(n) = n-1$ and $y(n) = n-2$.
Sequences satisfying this recurrence, but with different linear functions
  for $x(n)$ and $y(n)$ have been investigated by several authors
  in recent years, but the general behavior of these sequences remains
  rather mysterious.
Perhaps the most well-behaved sequences in the family occur when
  $x(n) = n$ and $y(n) = n-1$.
For a given parameter $s \ge 0$, we will show that the sequences 
  with $x(n) = y(n)+1 = n-s$ are almost as well-behaved.
In particular, we will show that they occur in a natural combinatorial
  setting, that they satisfy a recurrence relation of the
  form $a_s(n) = f(n) + a_s( n - g(n) )$, and that they have a
  fairly simple ordinary generating function.
  
The case of $s=1$ was studied before by Tanny \cite{Tanny}.
The case of $s=0$ was considered before by Conolly \cite{Conolly}.
Our attempt here is to simplify, unify, generalize, and
  combinatorialize their results.

\section{Meta-Fibonacci Sequences and Complete Binary Trees}

\begin{figure}
\includegraphics{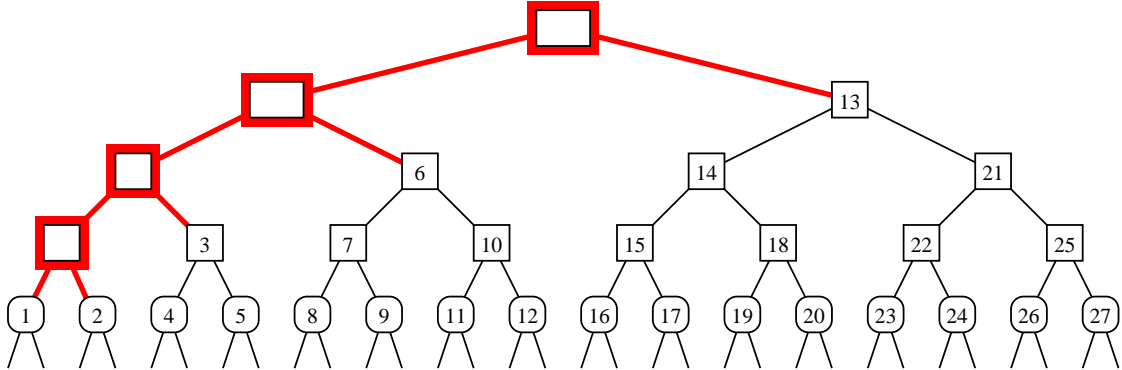}
\caption{The tree $\mathcal{F}_0$.}
\label{fig:F0}
\end{figure}

\begin{figure}
\includegraphics{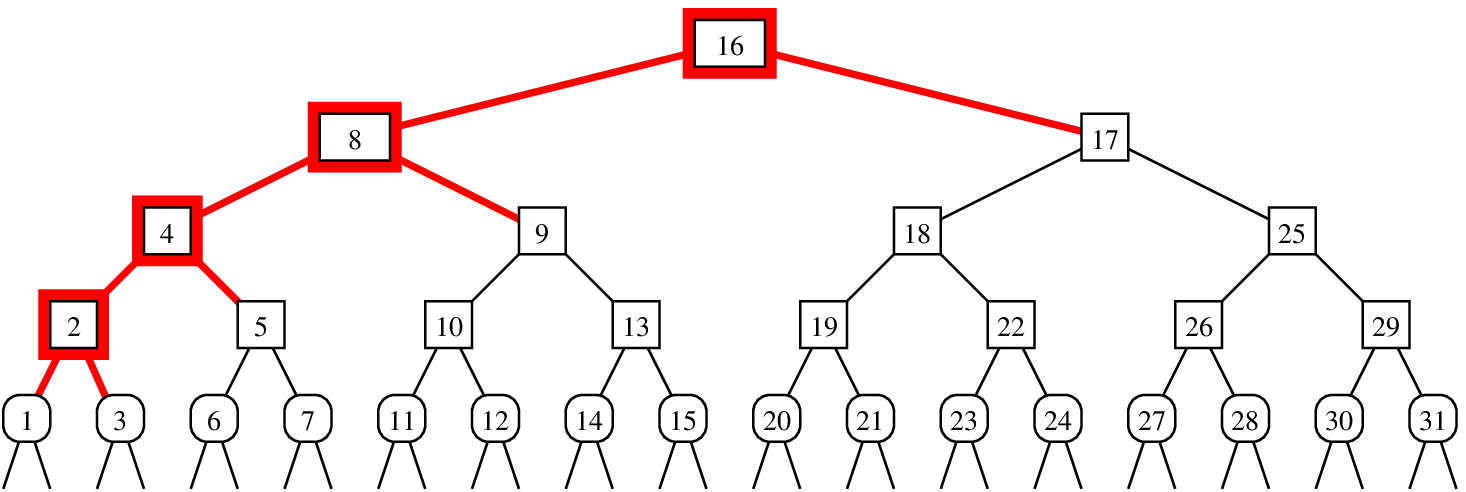}
\caption{The tree $\mathcal{F}_1$.}
\label{fig:F1}
\end{figure}

\begin{figure}
\includegraphics{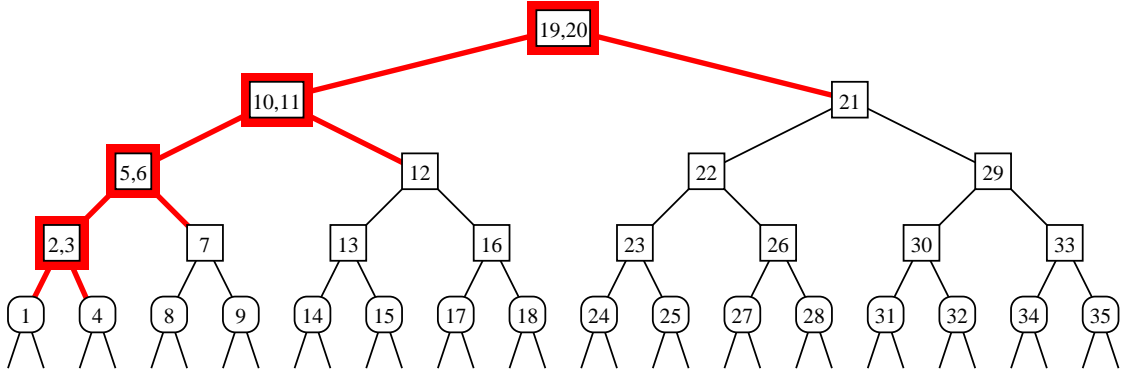}
\caption{The tree $\mathcal{F}_2$.}
\label{fig:F2}
\end{figure}

Figure \ref{fig:F0} shows part of an infinite ordered binary tree
  $\mathcal{F}$.
The forest of labelled trees in $\mathcal{F}$ consists of a succession 
  of complete binary
  trees of sizes $1,1,3,7,\ldots,2^h-1,\ldots$.
We refer to the subtree with $2^h-1$ nodes as \emph{subtree $h$}.
The nodes of these trees are labelled in preorder.
Now adjoin to $\mathcal{F}$ an infinite path that connects the trees
  from left-to-right as shown in Figure \ref{fig:F1}.
We will think of this path as being parameterized by a value $s$
  that gives the delay between the preorder counts of successive trees.
Alternatively, we can think of the nodes along the path as being
  super-nodes, where each super-node contains $s$ ordinary nodes.
This infinite tree is denoted $\mathcal{F}_s$, with our initial
  tree $\mathcal{F} = \mathcal{F}_0$. 
The trees $\mathcal{F}_0$, $\mathcal{F}_1$, $\mathcal{F}_2$ are
  shown in Figures \ref{fig:F0}, \ref{fig:F1}, \ref{fig:F2},
  respectively.  

Denote by $\T{s}(n)$ the tree induced by the first $n$
  nodes of the infinite tree $\mathcal{F}_s$.
Define $a_s(n)$ to be the number of nodes at the bottom level in 
  $\T{s}(n)$.
Also define $d_s(n)$ to be 1 if the $n$-th node is a leaf and to be
  0 if the $n$-th node is an internal node.
Finally, define $p_s(n)$ to be the positions occupied by the
  1's in the $d_s$ sequences.
Table \ref{table:numbers} gives the values of $a_s(n)$, $d_s(n)$,
  and $p_s(n)$ for
  $s = 0,1,2$ and $1 \le n \le 20$.
The values of four of these table entries appear in OEIS\footnote{OEIS =
  Neil Sloane's online encyclopedia of integer sequences.}, namely
  $a_0(n) = A046699$, $a_1(n) = A006949$, $d_0(n) = A079559$,
  and $p_0(n) = A101925 = A005187(n)+1$.
For fixed $s$ these numbers are related as follows.
\begin{equation}
a_s(n) = \sum_{j=0}^n d_s(n) \ \  \mbox{ and } \ \ 
p_s(n) = \min\{ j : a_s(j) = n \}.
\label{eq:sumd}
\end{equation}  

\begin{table}
\begin{tabular}{ccccccccccccccccccccccccccccc}
 $n = $&1&2&3&4&5&6&7&8&9&10&11&12&13&14&15&16&17&18&19&20 \\ \hline
$a_0(n)$ & 1&2&2&3&4&4&4&5&6&6&7&8&8&8&8&9&10&10&11&12 \\
$a_1(n)$ & 1&1&2&2&2&3&4&4&4&4&5&6&6&7&8&8&8&8&8&9 \\
$a_2(n)$ & 1&1&1&2&2&2&2&3&4&4&4&4&4&5&6&6&7&8&8&8 \\ \hline
$d_0(n)$ & 1&1&0&1&1&0&0&1&1&0&1&1&0&0&0&1&1&0&1&1  \\
$d_1(n)$ & 1&0&1&0&0&1&1&0&0&0&1&1&0&1&1&0&0&0&0&1  \\
$d_2(n)$ & 1&0&0&1&0&0&0&1&1&0&0&0&0&1&1&0&1&1&0&0  \\ \hline
$p_0(n)$ & 1&2&4&5& 8& 9&11&12&16&17&19&20&23&24&26&27&32&33&35&36 \\
$p_1(n)$ & 1&3&6&7&11&12&14&15&20&21&23&24&27&28&30&31&37&38&40&41 \\
$p_2(n)$ & 1&4&8&9&14&15&17&18&24&25&27&28&31&32&34&35&42&43&45&46 \\
\end{tabular}
\caption{The values of $a_s(n)$ and $d_s(n)$ for
  $s = 0,1,2$ and $1 \le n \le 20$.}
\label{table:numbers}
\end{table}

The $a_s(n)$ numbers satisfy the meta-Fibonacci recurrence relation
  stated in  Theorem 1 below.
  
\begin{theorem}
If $0 \le n \le s+1$, then $a_s(n) = 1$.
If $n = s+2$ then $a_s(n) = 2$.
If $n > s+2$, then
\[
a_s(n) = 
    a_s( n-s-a_s(n-1) ) + a_s( n-s-1-a_s(n-2) ).
\]
\end{theorem}
\begin{proof}
First observe that if all the leaves at the last level are
  removed from $\mathcal{F}_s$, then the same structure remains,
  except that the leftmost super-node needs to be made into an
  ordinary node (by subtracting $s-1$).
We will refer to this process as \emph{chopping} the last level.

We split the proof into two broad cases depending on whether
  $n$ is a leaf or not; i.e., whether $d_s(n) = 1$ (Case 1)
  or $d_s(n) = 0$ (Case 2).
  
\textbf{Case 1a:}
If $d_s(n-1) = d_s(n) = 1$ then $n$ and $n-1$ are sibling leaves
  and $a_s(n)$ is even.
For example, node 28 in Figure \ref{fig:F2}.
The trees $\T{s}(n-1)$ and $\T{s}(n-2)$ have the same number of
  nodes, $a_s(n)/2$, at the penultimate level as does $\T{s}(n)$.
Thus by chopping the last level from $\T{s}(n-1)$ and $\T{s}(n-2)$,
  we see that $a_s(n-s-a_s(n-1)) = a_s(n)/2 = a_s(n-s-1-a_s(n-1))$.


\textbf{Case 1b:}
If $d_s(n) = 1$ and $d_s(n-1) = 0$ then $n$ is a left child
  of its parent $n-1$ and $a_s(n)$ is odd.
For example, node 27 in Figure \ref{fig:F2}.
The tree $\T{s}(n-1)$ has $(a_s(n)+1)/2$ nodes at the penultimate
  level and the tree $\T{s}(n-2)$ has $(a_s(n)-1)/2$
  nodes at the penultimate level.
Thus by chopping the last level from $\T{s}(n-1)$ and $\T{s}(n-2)$,
  we see that $a_s(n-s-a_s(n-1)) = (a_s(n)+1)/2$ and
  $a_s(n-s-1-a_s(n-1)) = (a_s(n)-1)/2$.
  
\textbf{Case 2a:}
If $d_s(n) = 0$ and $d_s(n-1) = 1$, then $a_s(n)$ is even. 
For example, node 26 or node 29 in Figure \ref{fig:F2}.
The trees $\T{s}(n-1)$ and $\T{s}(n-2)$ have the same number of
  nodes, $a_s(n)/2$, at the penultimate level.
Node $n$ may have been at the penultimate level in $\T{s}(n)$, 
  but it is removed in $\T{s}(n-1)$ and $\T{s}(n-2)$.
Thus by chopping the last level from $\T{s}(n-1)$ and $\T{s}(n-2)$,
  we see that $a_s(n-s-a_s(n-1)) = a_s(n)/2 = a_s(n-s-1-a_s(n-1))$.

\textbf{Case 2b:}
If $d_s(n) = 0$ and $d_s(n-1) = 0$, then $a_s(n)$ is even.
For example, node 22 or node 30 in Figure \ref{fig:F2}.
The trees $\T{s}(n-1)$ and $\T{s}(n-2)$ have the same number of
  nodes, $a_s(n)/2$, at the penultimate level.
Node $n$ may have been at the penultimate level in $\T{s}(n)$, 
  but it is removed in $\T{s}(n-1)$ and $\T{s}(n-2)$.
Thus by chopping the last level from $\T{s}(n-1)$ and $\T{s}(n-2)$,
  we see that $a_s(n-s-a_s(n-1)) = a_s(n)/2 = a_s(n-s-1-a_s(n-1))$.
\end{proof}

Define $\mathcal{D}_s$ to be the infinite string $d_s(1) d_s(1) d_s(2) \cdots $.
Let $D_n$ be the finite string defined by $D_0 = 1$ and $D_{n+1} = 0 D_n D_n$.
Let $E_n$ be the finite string defined by $E_0 = 1$ and $E_{n+1} = E_n E_n 0$.

\begin{lemma}
\begin{equation}
\mathcal{D}_0 = D_0 D_0 D_1 D_2 D_3 \cdots = E_\infty
\label{eq:D0}
\end{equation}
\label{lemma:Einfty}
\end{lemma}
\begin{proof}
The first equality in (\ref{eq:D0}) is implied 
  immediately by the definition of $\mathcal{F}_0$;
  i.e., in $0 D_n D_n$ the 0 is from the root (which is listed first in preorder)
  and $D_n D_n$ is from the left and right subtrees.
By the definitions, $E_n^R = D_n$, where the superscript $R$ denotes reversal
  of the string.
Thus 
\[
D_0 D_0 D_1 \cdots D_n = E_n \cdots E_1 E_0 E_0.
\]
Since $E_n$ is a prefix of $E_{n+1}$ by definition, the expression
  $E_\infty$ is well-defined.
Hence $\mathcal{D}_0 = E_\infty$.
\end{proof}

The sequence $E_n$ has been considered before by 
  Allouche, Betrema, and Shallit \cite{ABS} in a different context.
It is interesting to note that the sequence $d(1) d(2) \cdots $
  is the limit of the morphism $0 \mapsto 0$ and $1 \mapsto 110$
  (also discussed in \cite{ABS}, pg. 237).
The following corollary is equation (6; pg. 132) in \cite{Conolly}.

\begin{corollary}
The numbers $a_0(n)$ satisfy the following recurrence relation
  for $0 \le k < 2^h$.
\[
a_0( 2^h-1+k ) = 2^{h-1} + a_0(k).
\]
\end{corollary}
\begin{proof}
Since $\mathcal{D}_0 = E_h E_h 0 \cdots$ and $|E_h| = 2^h-1$,
  the value of $d_0(2^h-1+k) = d_0(k)$ for
  $1 \le k \le 2^h-1$.
Since we defined $d_0(0) = 0$ it also holds when $k = 0$.
The number of 1's in $E_h$ is $\#_1(E_h) = 2^{h-1}$.
Thus
\begin{eqnarray*}
a_0(2^h-1+k) & = & \sum_{j=0}^{2^{h}-1} d_0(j) + \sum_{j=0}^k d_0(2^h-1+j) \\
& = & \#_1(E_h) + \sum_{j=0}^k d_0(j) \\
& = & 2^{h-1} + a_0(k).
\end{eqnarray*}
\end{proof}

\begin{lemma}
\[
a_s(n) = \begin{cases}
  a_0(n - sh) & \text{ if } 2^h + (s-1)h + 1 \le n \le 2^{h+1} + (s-1)h - 1, \\
  2^{h-1} & \text{ if } 2^h + (s-1)h - s + 1 \le n \le 2^h + (s-1)h.
\end{cases}
\]
\label{lemma:s0recurrence}
\end{lemma}
\begin{proof}
The labels on the nodes in subtree $h$ in $\mathcal{F}_s$ are
  exactly the values of $n$ lying in the first range above.
This is true since there are $1+1+3+\cdots+(2^{h-1}-1) = 
  2^h - h$ nodes in the subtrees to the left of subtree $h$,
  and $sh$ super-nodes.
Thus the lowest label of a node in subtree $h$ is
  $2^h-h+sh+1 = 2^h + (s-1)h + 1$, and the highest label
  is $2^h + (s-1)h + 2^{h}-1$.
The difference between the labels on corresponding nodes in
  $\mathcal{F}_s$ and $\mathcal{F}_0$ is $sh$ if the
  nodes are in subtree $h$; thus $a_s(n) = a_0(n-sh)$.

In the second range the nodes are super-nodes lying between
  subtrees $h-1$ and $h$ and therefore having $2^{h-1}$ leaves
  in their left-subtree. 
\end{proof}

\begin{corollary}
\[
a_1(n) = a_0( n - \lfloor \lg n \rfloor ).
\]
\end{corollary}
\begin{proof}
Taking $s = 1$ in Lemma \ref{lemma:s0recurrence} we obtain
  $a_1(n) = a_0(n-h)$ in the range $2^h+1 \le n \le 2^{h+1}-1$.
In that range $h = \lfloor \lg n \rfloor$.
We need only check what happens when $n = 2^h$.
By the lemma $a_1(2^h) = 2^{h-1}$.
However, in $\mathcal{F}_0$ the node $2^h-h$ is the rightmost
  node in subtree $h$ and thus $a_0(2^h-h) = 2^{h-1}$.
\end{proof}

The case $s = 1$ of the theorem below is roughly equivalent to
  equations (2.2) and (2.3) in Tanny \cite{Tanny}.
  
\begin{theorem}
If $1 \le k \le 2^{h-1}-1$, then
\[
a_s(2^h+(s-1)h+k+1) = 2^{h-2} + a_s(2^{h-1}+(s-1)h-s+k+1).
\] 
If $1 \le k \le 2^{h-1}-1$, then
\[
a_s(2^h+2^{h-1}+(s-1)h+k) = 2^{h-1} + a_s(2^{h-1}+(s-1)h-s+k+1).
\]
If $2^h+(s-1)h-s+1 \le n \le 2^h+(s-1)h+1$, then $a_s(n) = 2^{h-1}$. 
\end{theorem}
\begin{proof}
Let the node $n$ be in the subtree $h$ or the super-node, call it
  $y$, that is the parent of subtree $h$.
Let $x$ be the root of that subtree and denote the left and
  right subtrees of $x$ by $T_L$ and $T_R$.
We will prove the following recurrence relation. 
\begin{equation}
a_s(n) = \begin{cases}
2^{h-1} + a_s( n - 2^h - s + 1 ) & \text{ if $n \in T_R$}, \\
2^{h-2} + a_s( n - 2^{h-1} - s ) & \text{ if $n \in T_L$}, \\
2^{h-1} & \text{ if $n = x$ or $n \in y$}.
\end{cases}
\label{eq:treecases}
\end{equation}
Let $T$ be the subtree whose root is the right child of the
  left child of $y$.
In the first two cases above we are mappping the subtree $T_L$ or $T_R$
  to $T$, which has the same structure.
In the case of $T_R$ we skip over $2^{h-1}$ leaves and $2^h + s -1$
  nodes.
In the case of $T_L$ we skip over $2^{h-2}$ leaves and $2^{h-1}+s$ nodes.
Clearly, if $n = x$ or $n \in y$, then $a_s(n) = 2^{h-1}$.
  
From the proof of the previous lemma we know that $x = 2^h+(s-1)h+1$ and
  thus that the root of $T_R$ is $2^h + 2^{h-1}+(s-1)h+1$ and the
  root of $T_L$ is $x+1 = 2^h+(s-1)h+2$.
Thus we know the exact range of $n$ in each of the subtrees and the
  theorem statement is another way of writing (\ref{eq:treecases}).
\end{proof}
  
Let $r_1, r_2, r_3, r_4, \ldots = 1,2,1,3,1,2,1,4,1,2,1,3,1,2,1 \ldots$ 
  be the transition
  sequence of the binary reflected Gray code; this sequence is also known
  as the ``ruler function" (A001511).
If the alternating 0's are removed from the sequence
  $r_1-1, r_2-1, r_3-1, r_4-1, \ldots = 0,1,0,2,0,1,0,3,0,1,0,2,0,1,0 \ldots$
  then the ruler function is again obtained.
This implies that the generating function of the ruler function is
\begin{equation}
\sum_{k \ge 1} r_k z^k = \sum_{n \ge 0} \frac{z^{2^n}}{1-z^{2^n}}.
\label{eq:rulergf}
\end{equation}

\begin{lemma}
\begin{eqnarray*}
\mathcal{D}_0 & = & 11 0^{r_1} 11 0^{r_2} 11 0^{r_3} 11 0^{r_4} \cdots \\
              & = & 1 0^{r_1-1} 1 0^{r_2-1} 1 0^{r_3-1} 1 0^{r_4-1} \cdots 
\end{eqnarray*}
\label{lemma:DD0}
\end{lemma}
\begin{proof}
The ruler sequence is $R_\infty$ where $R_1 = 1$ and $R_{n+1} = R_n,n+1,R_n$.
Since $|R_n| = 2^n-1$, we have $r_{2^n+i} = r_i$ for $1 \le i \le 2^n-1$ 
  and $r_{2^n} = n+1$.
We will show that 
\[
E_{n+1} = 11 0^{r_1} 11 0^{r_2} \cdots 11 0^{r_{2^{n-1}}},
\]
which will finish the proof of the first equality since $\mathcal{D}_0 = E_\infty$.
By induction
\begin{eqnarray*}
E_{n+2} & = & E_{n+1} E_{n+1} 0 \\
& = & 11 0^{r_1} 11 0^{r_2} \cdots 11 0^{r_{2^{n-1}}} \ 
      11 0^{r_1} 11 0^{r_2} \cdots 11 0^{r_{2^{n-1}}} \  0 \\
& = & 11 0^{r_1} 11 0^{r_2} \cdots 11 0^{r_{2^{n-1}}} \ 
      11 0^{r_{1+2^{n-1}}} 11 0^{r_{2+2^{n-1}}} \cdots 110^{r_{2^{n}-1}} 11 0^{n+1},
\end{eqnarray*}
as required.
The second equality follows from the well-known property of the ruler
  sequence that $R_{\infty} = 1 + (0,r_1,0,r_2,0,r_3,0,r_4,0,\ldots )$.
\end{proof}

We can extend some of the previous results about $\mathcal{D}_0$ to
  $\mathcal{D}_s$.
For proposition $P$ the notation $\I{P}$ means 1 if $P$ is true
  and 0 if $P$ is false.
  
\begin{lemma}
Let $s_j = r_j + s \I{j \mbox{ is a power of 2}}$.
\[
\mathcal{D}_s = D_0 0^s D_0 0^s D_1 0^s D_2 0^s D_3 0^s \cdots
\label{eq:Ds}
\]
\[
\mathcal{D}_s = 1 0^{s_1-1} 1 0^{s_2-1} 1 0^{s_3-1} 1 0^{s_4-1} \cdots
\]
\label{lemma:Druler}
\end{lemma}
\begin{proof}
The proof is similar to those used in Lemmata \ref{lemma:DD0} 
  and \ref{lemma:Einfty} and is omitted.
\end{proof}

Since the $p_s(n)$ numbers give the positions of the 1's in 
  $\mathcal{D}_s$ the following corollary is true.
  
\begin{corollary}
For all $n \ge 1$,
\[
p_s(n+1) - p_s(n) = r_n + s \I{n \mbox{ is a power of 2}}.
\]
\label{corr:pdiff}
\end{corollary}  

\subsection{Generating Functions}

If $S = s(1)s(2) \cdots s(m)$ is a string then we use
  $S(z)$ to denote the ordinary generating function
  $S(z) = \sum s(i)z^i$.
Let $\mathcal{A}_s(z)$ and $\mathcal{D}_s(z)$ 
  denote the ordinary generating functions of the
  $a_s(n)$ and $d_s(n)$ sequences, respectively.
Directly from the definitions we get the equation shown below:
\[
\mathcal{A}_s(z) = \frac{\mathcal{D}_s(z)}{1-z}.
\]
Since $\mathcal{A}_s(z)$ is determined by 
  $\mathcal{D}_s(z)$ and $\mathcal{D}_s(z)$ is easier to
  treat, we first concentrate our attention on $\mathcal{D}_s(z)$.

\begin{lemma}
\[
D_n(z) = z^n (1+z)(1+z^3) \cdots (1+z^{2^n-1}) = z^n \prod_{j=1}^n (1+z^{2^j-1}).
\]
\[
E_n(z) = z (1+z)(1+z^3) \cdots (1+z^{2^n-1}) = z \prod_{j=1}^n (1+z^{2^j-1}).
\]\end{lemma}
\begin{proof}
From the recurrence relation $D_0 = 1$ and $D_{n+1} = 0 D_n D_n$ we obtain
  $D_0(z) = z$ and
\[
D_{n+1}(z) = z D_n(z) + z^{|0D_n|} D_n(z) = z(1+z^{2^{n+1}-1}) D_n(z).
\]
Similarly $E_0(z) = z$ and $E_{n+1}(z) = (1+z^{2^{n+1}-1}) E_n(z)$.
The result now follows by induction.
\end{proof}

\begin{corollary}
\[
\mathcal{D}_0(z) = z (1+z)(1+z^3)(1+z^7) \cdots  = z \prod_{n \ge 1} (1+z^{2^n-1}).
\]
\label{corr:D0(z)}
\end{corollary}
\begin{proof}
Follows at once from the the preceding lemma and the equation
  $\mathcal{D}_0 = E_\infty$ from Lemma \ref{lemma:Einfty}.
\end{proof}

\begin{theorem}
The generating function $\mathcal{D}_s(z)$ is equal to
\begin{equation}
z(1+z^{s+2^0}(1+z^{s+2^1}[1+z^{2^1-1}](1+z^{s+2^2}[1+z^{2^2-1}](1+z^{s+2^3}[1+z^{2^3-1}](1+  \cdots
\label{eq:Ds(z)}
\end{equation}
\end{theorem}
\begin{proof}
We need to translate the string
  $D_0 0^s D_0 0^s D_1 0^s D_2 0^s D_3 0^s \cdots$ from
  Lemma \ref{lemma:Druler} into its generating function.
Since 
\begin{equation}
|D_0 0^s D_0 0^s D_1 0^s \cdots D_{n-1} 0^s| = 
s+1 + \sum_{j=0}^{n-1} (2^{j+1} - 1 + s) = 2^{n+1} + (s-1)(n+1),
\label{eq:Dslength}
\end{equation}
we can write
\begin{equation}
\mathcal{D}_s(z) = z \left( 1 + \sum_{n \ge 0} z^{2^{n+1}+(n+1)(s-1)} D_n(z) \right).
\label{eq:Dsz}
\end{equation}
Let $x_k = z(1+z^{2^k-1})$, so that $D_n(z) = x_1 x_2 \cdots x_n$.
We can rewrite (\ref{eq:Ds(z)}) as
\begin{equation*}
z(1+z^{s+2^0}(1+z^{s+2^1-1}x_1(1+z^{s+2^2-1}x_2(1+z^{s+2^3-1}x_3(1+  \cdots
\end{equation*}
The coefficient of $x_1 x_2 \cdots x_n$ is $z$ raised to the power
  $1+ 2^{n+1} + (s-1)(n+1)$ by the sum given in (\ref{eq:Dslength}).
\end{proof}


\begin{theorem}
If $s \ge 1$, then
\begin{equation}
\mathcal{A}_s(z) = \frac{1-z^{s}}{1-z} \left( z + z \sum_{n \ge 1} 
  \prod_{k=1}^{n} z^{s-1} (z+z^{2^k}) \right)
\end{equation}
\label{thm:Dsz}
\end{theorem}

\begin{proof}
Call the expression on the right $R_s(z)$ and
  let $y = z^{s-1}$.
Multiply $R_s(z)$ by $1-z$, expand, and collect terms
  by increasing powers of $y$ to obtain

\begin{eqnarray*}
(1-z) R_s(z) & = & 
(1-zy) \left( z + z \sum_{n \ge 1} 
  \prod_{k=1}^{n-1} y (z+z^{2^k}) \right) \\
& = &
z + z \sum_{n \ge 1} 
  y^n \prod_{k=1}^{n} (z+z^{2^k}) 
- z^2 y - z^2 y \sum_{n \ge 1} 
  \prod_{k=1}^{n} y^n (z+z^{2^k})  \\
& = & 
z + z \sum_{n \ge 1} 
  y^n \left( \prod_{k=1}^{n} (z+z^{2^k}) 
- z y 
  \prod_{k=0}^{n} (z+z^{2^k}) \right) \\
& = & 
z + z \sum_{n \ge 1} 
  y^n \left( (z+z^{2^n}) \prod_{k=1}^{n-1}(z+z^{2^k}) 
- z 
  \prod_{k=0}^{n-1} (z+z^{2^k}) \right) \\
& = &   
z + z \sum_{n \ge 1} y^n z^{2^n}
  \prod_{k=1}^{n-1} (z+z^{2^k})
\end{eqnarray*}
Note that this last expression is equal to 
  $\mathcal{D}_s(z)$ by (\ref{eq:Dsz}).
\end{proof}

Jon Perry \cite{Perry} has observed experimentally that $a_1(n)$ counts the
  number of compositions of $n$ such that, for some $k$, 
\[
x_0 + x_1 + x_2 + \cdots + x_k = n
\ \  \mbox{ where } \ \ x_i \in \{1,2^i\}
\]
He uses the notation $1+[1,2]+[1,4]+[1,8]+ \cdots $
  to denote the set of such compositions and notes
  that many other combinatorial objects
  are in one-to-one correspondence with similar
  composition rules \cite{Perry}.
We call these rules \emph{specifications}.
  
\begin{corollary} For $s \ge 1$, the number of compositions of $n$ with
  specification
\[
[1,2,\ldots,s] + [s,2+s-1] + [s,4+s-1] + [s,8+s-1] + \cdots
\]
is $a_s(n)$.
\end{corollary}
\begin{proof}
This is clear from the generating function for $\mathcal{A}_s(z)$ given
  in Theorem \ref{thm:Dsz} once $z(1-z^s)/(1-z)$ is
  written as $z+z^2+\cdots+z^s$.
\end{proof}

As an example, for $s = 2$ and $n = 8$, the specification is
  $[1,2] + [2,3] + [2,5] + [2,9] + \cdots$ and the $a_2(8) = 3$ 
  compositions are
\[
8 = 1+2+5 = 1+3+2+2 = 2 + 2 + 2+ 2.
\]

To finish this section we also develop a 
  generating function for the $p_s(n)$ sequences.

\begin{lemma}
For all $s \ge 0$, 
\[
\sum_{n \ge 0} p_s(n) z^n =
\frac{1}{1-z} \left( 1 + z \sum_{k \ge 0} z^{2^k}\left( s + \frac{1}{1-z^{2^k}}
   \right) \right).
\]
\end{lemma}
\begin{proof}
Let $\mathcal{P}_s(z)$ denote the ordinary generating function of the
  numbers $p_s(n)$.
Then
\[
\sum_{n \ge 1} (p_s(n+1) - p_s(n)) z^n = \frac{1}{z}\left( (1-z) \mathcal{P}_s(z) -1 \right).
\]
By Corollary \ref{corr:pdiff} this expression is equal to 
\[
\sum_{n \ge 1} \left( r_n + s \I{n \mbox{ is a power of 2}} \right) z^n =
  \sum_{k \ge 0} \left( s z^{2^k} + \frac{z^{2^k}}{1-z^{2^k}} \right),
\]
where the equality follows from (\ref{eq:rulergf}).
Solving for $\mathcal{P}_s(z)$ finishes the proof.
\end{proof}

\section{Binary Compact Codes}

A binary compact code can be represented by an extended binary tree.  We use the term extended binary tree in the sense
  of Knuth \cite{Knuth1}: every node has either no children (a leaf) or two children (an internal node).  
Since no other types of codes are considered here, we shorten ``binary compact code" to ``code".  
A code of order $n$ can be represented by a tree with $n$ leaves in which the level numbers 
  $\ell_1 \ge \ell_2 \ge \cdots \ge \ell_n$ of the leaves are non-increasing.  
We will identify a compact binary code by the sequence of level numbers $(\ell_1,\ell_2,\ldots,\ell_n)$.  
For example, the codes for $n = 5$ are $(3,3,3,3,1)$, $(3,3,2,2,2)$, and $(4,4,3,2,1)$.  
Every code of order $n$ corresponds to a unique partition of 1 into the $n$
  powers of $1/2$ given by $1 = 2^{- \ell_1} + 2^{- \ell_2} + \cdots + 2^{- \ell_n}$.  
Thus $(3,3,3,3,1)$ corresponds to the partition $1 = 1/8 + 1/8 + 1/8 + 1/8 + 1/2$. 

The height $h$ of a tree is the length of the longest path from the root to any leaf.  
For a given height $h$ and integer $n$, we consider here the problem of finding the maximum number 
  of leaves at the largest level $h$ among all trees with $n$ leaves,
  which we denote by $M(n,h)$.  
Clearly $M(n,h) = 0$ if $h < \lceil \lg n \rceil$.  
A tree $T$ with $n$ vertices and height $h$ that has $M(n,h)$
  leaf pairs at the largest level is said to be an \emph{optimal tree}.  
We will show that our first two meta-Fibonacci sequences
  can be realized by certain families of optimal trees. 
This will be done via a ``greedy''
  algorithm for constructing a sequence of optimal tree/codes
  for successive values of $n$ and a fixed value $h$.
We denote these trees $\T{n,h}$ for natural numbers $n$ and $h$
  and call them \emph{greedy trees}.
Here is the greedy algorithm for constructing $\T{n,h}$.

\begin{itemize}
\item
If $n = h+1$, then there is only one tree/code,
  namely $h,h,h-1,\ldots,2,1$.
\item
Given $\T{n,h}$ the code $\T{n,h}$ is obtained by replacing
  the leftmost level $\ell_i$ for which $\ell_i < h$ by the
  two levels $\ell_i+1,\ell_i+1$.
\end{itemize}

We will also consider the trees $\T{n} = \T{n,\lceil \lg n \rceil}$.
They may be constructed greedily as follows.
 
\begin{itemize}
\item
If $n=0$, then the tree is a leaf.  
\item
If $n = 2^h$, then $\T{n}$ is a complete binary tree (all
  leaves are at level $h$).
Tree $\T{n+1}$ is the tree of 
  height $h+1$ whose left subtree is $\T{n}$ and whose right
  subtree is a single leaf.
\item
If $n$ is not a power of 2, then expand the leftmost leaf which is
  not at the largest level, as described above.
\end{itemize}

It is also interesting to consider the inverse process of obtaining
  $\T{n}$ from $\T{n+1}$.
The inverse rule is very simple: 
  Replace the rightmost equal pair $\ell_j = \ell_{j+1}$ by 
  $\ell_j-1$.  
  
We could also have defined a code by the number of internal
  nodes at each level in the corresponding tree.
Given a code of height $h$, let $[\tau_0,\tau_1,\ldots,\tau_{h-1}]$
  be the sequence in which $\tau_i$ is the number of internal
  nodes at level $h$.
For our example codes given earlier, the corresponding level counts
  are $[1,1,2], [1,2,1], [1,1,1,1]$.
These counts clearly must satisfy
\[
\tau_i \le 2 \tau_{i-1}, \mbox{ for } 1 \le i \le h-1 \mbox{ and }
\]
\[
\tau_0 + \tau_1 + \cdots + \tau_{h-1} = n-1.
\]
Subject to these two constraints $M(n,h)$ is the largest
value that $\tau_{h-1}$ can attain.

Let $k$ be the largest level for which $\tau_k < 2 \tau_{k-1}$.
The greedy algorithm simply replaces $\tau_k$ by $1 + \tau_k$.

\begin{lemma}
Let $\tau_0,\tau_1, ... ,\tau_{h-1}$ be the number of vertices at each level for the
tree $\T{n,h}$ and suppose 
that $t_0,t_1, ... ,t_{h-1}$ are the vertex numbers by level for any other 
tree $T$ with $n$ leaves and height $h$. 
For any $0 \le j \le h-1$, 
\[
\tau_j + \cdots + \tau_{h-1} \ge t_j + \cdots + t_{h-1}.
\]
\end{lemma}
\begin{proof}
For any $h$, the result is true for $n = h + 1$ since there 
is only one tree with $h+1$ leaves and height $h$. 
Similarly it is true for $n = 2^h$.

Assume the result is true for all trees with $n$ (for some $n < 2^h$) leaves and height $h$.  
Let $\tau'_1,\tau'_2, \ldots ,\tau'_{h-1}$ be the vertex numbers by level for $\T{n+1,h}$.  
Note that for some $k > 0$, we have $\tau'_k = 1 + \tau_k$ and $\tau'_j = \tau_j$ for all $j \neq k$. 
By the greedy algorithm $\tau_{h-1} = 2 \tau_{h-2}, \ldots , \tau_{k+1} = 2 \tau_k$, 
  but $\tau_k < 2 \tau_{k-1}$. 
Let $T'$ be any tree with $n+1$ leaves and height $h$ and suppose 
  that $T$ is the tree with $n$ leaves and height $h$ formed by removing from $T'$ the 
  rightmost pair of leaves at level $h$. 
Suppose that $t_0,t_1, \ldots ,t_{h-1}$ are the vertex numbers by level for $T$.  
By induction, we assume that 
  $t_j + \cdots + t_{h-2} + t_{h-1} \le \tau_j + \cdots + \tau_{h-2} + \tau_{h-1}$
  for all $0 \le j \le h-1$. 
Let $t'_0,t'_1, \ldots ,t'_{h-1}$ be the vertex numbers by level for for $T'$. 
Note that $t'_{h-1} = 1 + t_{h-1}$ and $t'_j = t_j$ for all $j < h-1$. 
For all $i \le k$, we see that $t'_i + \cdots + t'_{h-1} \le \tau'_i + \cdots + \tau'_{h-1}$. 

Suppose $t'_i + \cdots + t'_{h-1} > \tau'_i + \cdots + \tau'_{h-1}$ for some $i > k$. 
Let $m \ge i$ be the smallest index such that $t'_m > \tau'_m = \tau_m$.
Since $m > k$ we have $\tau_m = 2 \tau_{m-1} = 2 \tau'_{m-1}$.
Since $m$ was chosen to be smallest $\tau'_{m-1} \ge t'_{m-1}$. 
Putting these inequalities together we have 
  $t'_m > 2 t'_{m-1}$ which is a contradiction. 
Thus $t'_i + \cdots + t'_{h-1} \le \tau'_i + \cdots + \tau'_{h-1}$ for all $i$ and the result 
  is true by induction.
\end{proof}
  
\begin{theorem}
The greedy algorithm produces optimal trees.
\end{theorem}
\begin{proof}
Let $t_{h-1}$ be the number of internal nodes in some code with $n$ leaves and height $h$.
The previous lemma tells us that $\tau_{h-1} \ge t_{h-1}$, where $\tau_{h-1}$ is the number
of internal nodes at level $h-1$ in $\T{n,h}$.
\end{proof}

\begin{corollary}
The largest number of 1's in a partition of $2^h$ into powers of 2 consisting of $n$ parts is $M(n,h)$.
\end{corollary}
\begin{proof}
Multiply the partition of 1 described above by $2^h$ to obtain a partition of $2^h$.
\end{proof}  

Define $a(n)$ to be the maximum number of leaf pairs at the largest level, 
  taken over all binary trees with $n$ leaves.  
In other words, $a(n) = \max{ M(n,h)\ :\ 1 \le h \le n-1}$.  

\begin{corollary}
For any $n$, we have $a(n) = M( n, \lceil \lg n \rceil ) = a_1 (n-1)$.
\end{corollary}  
\begin{proof}
For any $k > h = \lceil \lg n \rceil$ construct a optimal tree of height $k$ with $n$ 
  leaves using the greedy algorithm. 
The tree $\T{n,k}$ has a subtree attached to an interior vertex at level $k-h$ 
  which is ismorphic to $\T{n-k+h,h}$.  
Clearly $M(n-k+h,h) \le M(n,h)$ and thus $a(n) = M( n, \lceil \lg n \rceil )$.  
The tree $\mathcal{T}_1(n)$ from Section 2 which defines the sequence $a_1(n)$ has $n+1$ leaves 
  (since it has $n$ interior vertices) and is equal to the greedy tree 
  for $M(n+1,\lceil \lg (n+1) \rceil)$.  
Only the order in which the vertices are added is different since we add the vertices from 
  the bottom in constructing $\mathcal{T}_1(n)$.
\end{proof}

Similarly define $b(n)$ by the equation $b(n) = M(n+h,h)$ for $h+1 \le n+h \le 2^h$.  

\begin{corollary}
The sequence $b$ is well-defined and $b(n) = a_0 (n)$.
\end{corollary}  
\begin{proof}
For a given $n$, let $h$ be the smallest height such that $n+h \le 2^h$.  
For any larger height $k > h$ construct a optimal tree of height $k$ with $n+k$ leaves 
  using the greedy algorithm.  
The tree $\T{n+k,k}$ has a subtree attached to an interior vertex at level $k-h$ 
  which is equal to $\T{n+h,h}$.  Thus $M(n+k,k) = M(n+h,h)$.  
If $h$ is the height of the $n$th subtree $\T{n}$ which defines the sequence $a_0$, 
  we see inductively that $\T{n}$ has $n+h-1$ internal nodes and thus $n+h$ leaves.  
As before we see that $T_n$ is equal to the greedy tree for $M(n+h,h)$. 
Thus $\T{n}$ has $b(n)$ leaf pairs at the largest level $h$, and we see that $b(n) = a_0 (n)$.
\end{proof}  

Thus we have shown that the first 2 meta-Fibonacci sequences in our family 
  of sequences have concrete realizations as the solutions of 
  optimization problems involving binary compact codes/trees.  
        

\section{Acknowledgements}

We wish to thank Don Knuth, Jon Perry, Jeff Shallit, and Herb Wilf for helpful
  comments related to this research.

\end{document}